\documentstyle[12pt,amssymb]{article}
\begin{document}

\def\bR{\ifmmode{\Bbb{R}}\else{{$\Bbb{R}$ }}\fi}
\def\bZ{\ifmmode{\Bbb{Z}}\else{{$\Bbb{Z}$ }}\fi}
\def\bQ{\ifmmode{\Bbb{Q}}\else{{$\Bbb{Q}$ }}\fi}
\def\bN{\ifmmode{\Bbb{N}}\else{{$\Bbb{N}$ }}\fi}
\def\bC{\ifmmode{\Bbb{C}}\else{{$\Bbb{C}$ }}\fi}

\renewcommand{\theequation}{\thesection.\arabic{equation}}

\renewcommand{\thefootnote}{\fnsymbol{footnote}}

\rightline{FTUV 95-1}
\rightline{IFIC 95-1}

\begin{center}
{\Large{\bf An introduction to quantum groups and  non-commutative differential
calculus\footnote{Invited lecture at the {\it III Workshop on Differential
Geometry}, Granada (September 1994)}
}}
\end{center}

\vspace{1\baselineskip}

\begin{center}
{\large{J.A. de Azc\'{a}rraga$^{\dagger}$ and
F. Rodenas$^{\dagger \; \ddagger}$}}
\end{center}

\noindent
{\small{\it $ \dagger \;$ Departamento de F\'{\i}sica Te\'{o}rica and IFIC,
Centro Mixto Universidad de Valencia-CSIC,
E-46100-Burjassot (Valencia), Spain.}}

\noindent
{\small{\it $ \ddagger \;$ Departamento de Matem\'atica Aplicada,
 Universidad Polit\'ecnica de Valencia,
E-46071 Valencia, Spain.}}

\section{Introduction}

$\quad$ {\it Quantum groups} provide an interesting  example of non-commutative
geometry  \cite{CONNES} (for a review with theoretical physics in mind
see, {\it e.g.}, \cite{COQUE}). They
may be looked at in various ways \cite{COL}.  From a mathematical point of
view, they may be introduced by making emphasis on their
$q$-deformed enveloping algebra aspects \cite{DRINFELD,JIMBO}, which leads to
the {\it quantized universal enveloping algebras}, or by making
emphasis in the $R$-matrix formalism that describes the deformed group algebra
\cite{FRT}. A useful point of view for possible
physical applications
is to look at quantum groups as a generalization (deformation)
of symmetry groups
acting on  generalized representation spaces or {\it quantum spaces}
\cite{MANIN} (see also \cite{FRT,WZ,DAVID}). Quantum groups are
mathematically
well defined in the framework of Hopf algebras \cite{ABE};
their algebraic properties
depend on one deformation parameter $q$ (or more) in such a way that  for $q$=1
the deformed structures become the  standard (non-deformed, Lie)  ones.
Thus, the essential feature in the field of   quantum groups (we shall not
discuss  their dual {\it quantum algebra} aspect above)
is in some sense similar to the relation between classical
and quantum mechanics, where the commutative algebra of functions on phase
space (the algebra of observables) becomes non-commutative after quantization.
In the case of Lie groups, the commutative algebra of functions on the
group manifold is replaced by a non-commutative algebra after
quantization (or $q$-deformation);
in particular, the matrix elements generating the  algebra become
non-commutative.
This analogy justifies the  `quantum' name given to these structures.
It must be said, nevertheless, that
in physics `quantum' refers to the appearance of the Planck
constant  and that its relation to
$q$  is, at best, unclear. Even the analysis
of the `quasiclassical limit', usually introduced by assuming that
$q=e^{\hbar}$, requires writing $q=e^{\gamma \hbar}$ {\it if} $\hbar$ is
the Planck constant, {\it i.e.}, it introduces a new dimensional constant
$\gamma$ in classical physics.

\section{Quantum groups as $q$-symmetries}

\subsection{$GL_q$ and  quantum matrices}

$\quad$ Let us start by writing down some simple algebraic
quantum group aspects relevant for our discussion.
An element of $GL(2, \bC )$ is a   regular 2$\times$2 matrix
\begin{equation}\label{cero}
T = \left( \begin{array}{cc}
a & b\\
c & d
\end{array} \right)
\end{equation}
\noindent
with $a,b,c,d \in \bC$.
The {\it quantum group} $GL_{q} (2)$ is the associative
{\it algebra} (quantum groups are not really group manifolds)
generated by the entries  $a,b,c,d$ of a matrix $T$ satisfying the homogeneous
quadratic relations ($\lambda \equiv q-q^{-1}$, $q \in \bC$, $q \neq 0$)

\begin{equation}\label{ua}
\begin{array}{lll}
ab = qba\;, \quad & bd = qdb \;,\quad & bc = cb \;,\\
ac = qca \;,\quad & cd = qdc \;, \quad & [a,d] = \lambda bc \;.
\end{array}
\end{equation}
\noindent
By this reason, the notation $Fun(GL_q(2))$ is also used.
The relations (\ref{ua}) are not arbitrary; they are the result of certain
requirements. The most important properties of the matrices $T$
satisfying (\ref{ua}),
or {\it quantum matrices}, are (see \cite{VZW,DAVID}):

\,
\noindent
1) The matrix (co)multiplication  preserves (\ref{ua}).
If we take a second matrix
$T'$  with   entries satisfying  (\ref{ua}) and commuting  with  those of $T$,
the entries $T''=TT'$ satisfy (\ref{ua}) again.
In contrast, the entries of $T^n$
satisfy  (\ref{ua}) after replacing $q$ by $q^n$
(quantum matrices do not form a group);
this product should not be confused with the comultiplication $T''=TT'$ above.

\,
\noindent
2) The relations (\ref{ua}) are consistent, {\it i.e.}, they do not
generate higher order relations. For example, let us take the
same relations as in
(\ref{ua}) except for $bc=qcb$. If we put $acd$ into the order
$dca$ by exchanging $ac$  first  or by moving first $cd$ and compare a new,
third order relation (if $q  \neq 1$) $bc^2=0$ appears.
The form of the equations in  (\ref{ua}) guarantees that no such
higher order relations arise; they lead to a finitely generated quadratic
algebra.

 To go from $GL_q(2)$ (eqs. (\ref{ua})) to $SL_q(2)$ we have to remove
one generator. The element (see (\ref{adet}))

\begin{equation}\label{ub}
ad - qbc = da - q^{-1} bc : = det_{q} T \quad,
\end{equation}

\noindent
is a central (commuting) element of the algebra
which defines the $q$-determinant
of the matrix  $T$; the addition of the
constraint $det_{q} T = 1$ to eqs.
(\ref{ua})  consistently reduces the number of generators to three.
In the (`classical') limit $q$$=$$1$, (\ref{ua}) just expresses that the
algebra generated
by the elements of $SL(2, \bC)$ is commutative, and (\ref{ub}) is the usual
determinant.

The above discussion does not make  apparent why eqs. (\ref{ua})  plus
$det_{q} T$=$1$ define $SL_{q} (2)$, nor how to generalize
them to the $SL_{q} (n)$ case. The $q$-group {\it structure}
becomes clearer by using the $R$-matrix formalism originally developed in the
framework  of the quantum inverse scattering method (see \cite{FRT} and
references therein).
The non-commutativity of the entries of $T$ may be expressed by saying that
\begin{equation}
T_1T_2= \left[ \begin{array}{cccc}
                 aa & ab & ba & bb \\
                 ac & ad & bc & bd \\
                 ca & cb & da & db \\
                 cc & cd & dc &dd
                 \end{array} \right]   \neq  \left[ \begin{array}{cccc}
                 aa & ba & ab & bb \\
                 ca & da & cb & db \\
                 ac & bc & ad & bd \\
                 cc & dc & cd &dd
                 \end{array} \right]  = T_2T_1
\end{equation}
\noindent
where $T_{1} = T \otimes I$, $T_{2} = I \otimes T$ (see Appendix for
notation). Then,  eqs. (\ref{ua}) may be rewritten as `RTT'
(or `FRT' \cite{FRT}) relations,

\begin{equation}\label{uc}
R_{12} T_{1} T_{2} = T_{2} T_{1} R_{12} \quad , \quad
(R_{ij,ab} T_{ak} T_{bl} = T_{jc} T_{id} R_{dc,kl}) \quad,
\end{equation}
\noindent
where $R_{12}$ is the 4$\times$4  numerical matrix given in (\ref{af2}).
In this form, they may
be generalized to any dimension; all is needed is the
appropriate $n^{2} \times
n^{2}$ $\, R$-matrix which  for $GL_q(n)$ is \cite{FRT}
\begin{equation}\label{1j}
R_{ij,kl} = \delta_{ik} \delta_{jl} (1 + \delta_{ij} (q-1)) + \lambda
\delta_{il} \delta_{jk} \theta (i-j) \quad i,j...=1...n
\end{equation}
$$
\theta (i-j) = \left\{ \begin{array}{l}
                      0 \quad i \leq j \\
                      1  \quad i > j
                      \end{array} \right.  \quad.
$$

\noindent
The central $GL_q(n)$  $q$-determinant of $T=(t_{ij})$ ($i,j=1,...,n$) is given
\cite{FRT} by
\begin{equation}
det_qT= \sum_{s \in S_n} (-q)^{l(s)}\, t_{1\,s(1)}  ... t_{n\,s(n)}
 \quad ,
\end{equation}
\noindent
where  $l(s)$ is the `parity' of the permutation $s$;
the  condition $det_qT=1$ defines $SL_q(n)$.
Moreover, there exist matrices $R$ which define  the $q$-deformation of all
$A_l,B_l,C_l,D_l$-type simple groups and their corresponding real forms
as well as for the exceptional groups and supergroups; we refer to
\cite{FRT} for details.

One could insert other  matrix as $R$ in (\ref{uc}) in order to reproduce
(\ref{ua}).
However, the natural  ones (as (\ref{af2})) satisfy
the {\it Yang-Baxter equation} (YBE)
\begin{equation}\label{YBE}
R_{12}R_{13}R_{23}=R_{23}R_{13}R_{12}
\end{equation}
\noindent
which ensures the consistency of (\ref{uc}) (for other aspects of the YBE see
\cite{YB,COL} and references therein).
This means that no further relations for the
generators higher than the quadratic ones  (\ref{ua})
may be derived from (\ref{uc}) and the requirement of associativity
for the algebra, which  is postulated from the very beginning
and is independent of (\ref{YBE}). This equation
is sometimes introduced by reordering $T_1T_2T_3$ to
$T_3T_2T_1$ by two different paths using the RTT  relation and the
associativity property of the algebra. In this way
one is lead to $(R_{12}R_{13}R_{23}-R_{23}R_{13}R_{12})T_1T_2T_3=
T_3T_2T_1(R_{12}R_{13}R_{23}-R_{23}R_{13}R_{12})$. Thus, eq. (\ref{YBE})
is {\it consistent} with eq. ({\ref{uc}), but it is not implied by it. To
see this explicitly, consider eq. (\ref{uc}) rewritten in the form
$\hat{R}T_1T_2=T_1T_2 \hat{R}$ using (\ref{af}) ($\hat{R}={\cal P}R$,
$\hat{R}_{ij,kl}=R_{ji,kl}$). Then, due to (\ref{af3})
we get ($q^2+1 \neq 0$)
\begin{equation}\label{proj}
P_+T_1T_2P_-=0 \quad , \quad P_-T_1T_2P_+=0 \quad.
\end{equation}
\noindent
The first equation implies $ab-qba$=0, $cd-qdc$=0 and $[a,d]$=$qbc-q^{-1}cb$,
while the second gives $ac-qca$=0, $bd-qdb$=0 and $[a,d]$=$qcb-q^{-1}bc$.
In all, these equations reproduce (\ref{ua}). These equations also follow
from $P_{ \pm}T_1T_2$=$ T_1T_2P_{\pm}$, {\it i.e.} from an `RTT'
relation with $P_{\pm}$
as an $\hat{R}$-matrix,  although ${\cal P}P_{\pm}$ are not invertible
and does not satisfy the YBE (\ref{YBE}).

Since quantum groups are very close to the algebra of functions on a Lie
group, we may expect them to have other characteristics pertaining to the
group multiplication rule, inverse (antipode) and unit elements, etc. In fact,
 they may be characterized as Hopf algebras
\cite{DRINFELD,JIMBO,FRT} (see \cite{TAK}
for  reviews), but this aspect will not be considered here.

\subsection{The quantum plane}

$\quad$ Let us now  introduce a
deformed `representation space' for $GL_q(2)$ (and hence for $SL_q(2)$).
This is the {\it quantum plane} $C_{q}^{2}$, or associative
{\it algebra} (a $q$-plane is not a manifold) generated by two elements
$(x,y)=X$ (a `$q$-two-vector') subjected to the commutation property
\cite{MANIN}
\begin{equation}\label{1a}
xy = qyx \quad .
\end{equation}
\noindent
The commutation relation (\ref{1a}) can also be expressed  by using the
$q$-symplectic metric $\epsilon^q$ \cite{VZW,DAVID}
\begin{equation}\label{epsilon1}
\epsilon^q = \left( \begin{array}{cc}
                       0 & q^{-1/2} \\
                       -q^{1/2} & 0
                     \end{array} \right)  \quad , \qquad
( \epsilon^q )^2=-I
\end{equation}
\noindent
by the  equation
\begin{equation}\label{epsilon2}
X^t \epsilon^q X =0 \quad,\qquad \epsilon^q_{ij} X_iX_j=0
\end{equation}
\noindent
which reflects that the $q$-symplectic norm of a $q$-two-vector  vanishes.
It is also possible to introduce a pair of (odd) variables $(\xi,\eta)=\Omega$
(an odd $q$-two-vector)  satisfying
\begin{equation}\label{odd}
\xi \eta = - \frac{1}{q} \eta \xi \quad ,\qquad \xi^2 =0= \eta^2 \quad;
\end{equation}
\noindent
for $q$=1, $(x,y)$ commute and $(\xi,\eta)$ anticommute
(in a non-commutative differential calculus this second set of variables
may be  identified \cite{WZ,MANIN2} with the differentials  of $(x,y)$).
When  it is required that  after the transformation
({\it coaction}\footnote{Specifically, the
mapping (coaction) $\varphi : C^2_q \mapsto GL_q(2) \otimes C^2_q$,
$\varphi(X_i)= T_{ij}  \otimes X_j$   is an algebra
homomorphism, and $C^2_q$  is a left $GL_q(2)$-{\it comodule}, see below.})
$X'=TX$, $\Omega'=T \Omega$ (the entries
of $T$ commute with those of  $X$ and $\Omega$) the new entities
$(x',y')$, $(\xi',\eta')$ satisfy also (\ref{1a}), (\ref{odd}),
the commutation properties of the elements of $T$
are completely determined and (\ref{ua}) is obtained.
This allows us to consider the quantum plane (\ref{1a})
as the `representation' space
of the $GL_q(2)$ quantum group  (\ref{ua}).
Since  the
non-commuting properties of the quantum group are encoded in the $R$-matrix
by (\ref{uc}), it is natural  to
define the non-commuting properties of the $q$-plane analogously. Indeed,
eqs. (\ref{1a}) and (\ref{odd})  may be expressed as
\begin{equation}\label{1d}
R_{12} X_{1} X_{2} =  q X_{2} X_{1}
\quad \Longleftrightarrow \quad
R_{21}^{-1} X_{1} X_{2} =  q^{-1} X_{2} X_{1}\quad,
\end{equation}
\begin{equation}\label{1d1}
\quad R_{12} \Omega_{1} \Omega_{2} =  -q^{-1} \Omega_{2} \Omega_{1}
\quad \Longleftrightarrow \quad
R_{21}^{-1} \Omega_{1} \Omega_{2} =  -q \Omega_{2} \Omega_{1} \quad ,
\end{equation}
\noindent
where $X_{1} X_{2}$ and
$X_{2} X_{1}$ are, respectively, the four-vectors ($xx,xy,$ $yx,yy$) and
($xx,yx$, $xy,yy$) (analogously for $\Omega_{1} \Omega_{2}$ and
$\Omega_{2} \Omega_{1}$). In components, (\ref{1d})  reads
$\,R_{ij,kl} X_{k} X_{l} = q X_{j} X_{i}$, $\,\hat{R}_{ij,kl} X_{k} X_{l}
= q X_{i} X_{j}$;  similar expressions  are obtained for eq. (\ref{1d1}).
Both relations in (\ref{1d}) (and in ({\ref{1d1}))  are equivalent
since  $({\cal P}R{\cal P})_{ij,kl}=R_{ji,lk}$
({\it i.e.}, ${\cal P}R_{12}{\cal P}=R_{21})$  and $({\cal
P}X_1X_2)_{ij}=(X_1X_2)_{ji}$
{\it i.e.}, ${\cal P}X_1X_2=X_2X_1$.
Eqs. (\ref{1d}) and (\ref{1d1}) are preserved by the $q$-transformations
$X'=TX$ and $\Omega'=T \Omega $
since the components of $X$ and $\Omega$ are assumed to
commute with the entries of $T$. For instance,
\begin{equation}\label{1f}
\begin{array}{l}
R_{12} X_{1}' X_{2}' = R_{12} (T_{1} X_{1})
(T_{2} X_{2})=
R_{12} T_{1} T_{2} X_{1} X_{2} \\
= T_{2} T_{1} R_{12}
X_{1} X_{2} =q T_{2} T_{1} X_{2} X_{1} = q X_{2}' X_{1}'\quad
\end{array}
\end{equation}

\noindent
using (\ref{uc}): the invariance of the commutation properties  (\ref{1d})
under a  `$q$-symmetry' transformation requires (\ref{uc}).
Although the indices in all previous expressions take the values $1, 2$,
the $R$-matrix form of
(\ref{1d}) and (\ref{uc}) makes it clear how to
generalize them to $GL_{q}(n)$; all that is needed is the appropriate
$n^{2} \times n^{2}$  $R$-matrix  given by eq. (\ref{1j}).
With it, the relations defining the `quantum hyperplane'
\begin{equation}\label{1k}
X=(x_{1},...,x_{n})\quad , \quad
x_{i} x_{j} = q x_{j} x_{i} \quad (i<j)\; i,j=1...n
\end{equation}

\noindent
are again expressed by (\ref{1d}) and preserved under $GL_{q}(n)$
because of (\ref{uc}).

All this is  well known. Let us now  show
how to extend  these $q$-vector constructions
(see also \cite{PRAGA} and references therein).
We shall consider here the simplest example of $q$-twistors
constructed from $q$-two-vectors (spinors) [(\ref{1a}),
(\ref{epsilon2}), (\ref{1d})] and
the  application to $q$-Minkowski space algebras.

\setcounter{equation}{0}

\section{Beyond $q$-vectors}

\subsection{Second rank $q$-tensors. $q$-twistors}

$\quad$ Since a two-dimensional object transforming under $SL(2,\bC)$ is a
spinor, an element of the $q$-plane should be called a $q$-{\it spinor}
rather than $q$-vector. Since the complex conjugation of $A \in SL(2, \bC)$
defines an
unequivalent representation, spinors come in two varieties, dotted and
undotted. The classical construction of a (Minkowski) real four-vector uses
both,
\begin{equation}\label{3a}
K_{\alpha \dot{\beta}} =
(\sigma_{\mu} x^{\mu})_{\alpha \dot{\beta}}\quad \alpha, \dot{\beta}=1,2 \quad,
\quad K= \sigma_0 x^0+ \sigma_i x^i= \left( \begin{array}{cc}
                  x^0 + x^3 & x^1-ix^2 \\
                  x^1+ix^2  & x^0-x^3
                \end{array} \right) \quad,
\end{equation}

\noindent
and $K'_{\alpha \dot{\beta}} = A_{\alpha} \, .^{\gamma} K_{\gamma \dot{\delta}}
(\tilde{A}^{-1})^{\dot{\delta}} . _{\dot{\beta}} \; ,$ where $A$ and
$\tilde{A} \equiv  (A^{-1})^{\dagger}$ are the two fundamental
representations of
$SL(2, \bC)$, $D^{\frac{1}{2},0}$ and $D^{0,\frac{1}{2}}$; $\sigma^0=I_2$,
$\sigma^i$ ($i$=1,2,3) are the Pauli matrices and
$detK=(x^0)^2- \vec{x}^2=detK'$ is the square of Minkowski length.  Thus, a
$q$-deformation of the Lorentz group, $L_q$, may be obtained
\cite{POWO}-\cite{OSWZ}
by replacing the classical ($q$=1)
$A$ and $\tilde{A}$ by two copies $T$ and $\tilde{T}$ of
$SL_{q} (2)$ plus the reality condition ($*$-structure)
$T^{\dagger}=\tilde{T}^{-1}$.
The commutation relations in this general situation may be expressed in terms
of   four $R$-matrices $R^{(i)}$, $i=1,...,4$,
\begin{equation}\label{2b}
\begin{array}{ll}
R^{(1)} T_{1} T_{2} =  T_{2} T_{1} R^{(1)}\quad , &
\quad T_{1}^{\dagger} R^{(2)} T_{2} =  T_{2} R^{(2)} T_{1}^{\dagger}\quad ,\\
R^{(4)} T_{1}^{\dagger} T_{2}^{\dagger} =  T_{2}^{\dagger} T_{1}^{\dagger}
R^{(4)}\quad , & \quad
T_{2}^{\dagger} R^{(3)} T_{1} =  T_{1} R^{(3)} T_{2}^{\dagger}\quad,
\end{array}
\end{equation}
\noindent
where $R^{(2)}=R^{(3) \,\dagger}$, $R^{(1)}$ is the $SL_q(2)$ $R$-matrix
($R^{(1)}$ may be taken as  $R_{12}$ or $R_{21}^{-1}$) and,
since $\tilde{T}= (T^{\dagger})^{-1}$
is another copy of $SL_q(2)$, consistency requires
$R^{(4)}={\cal P} R^{(1)} {\cal P}=R^{(1)\,t}$  and $q$ real (from now on,
we shall take $q \in \bR$).
The matrix $R^{(2)}$ (and hence $R^{(3)}$) defines how the elements of both
quantum groups $T$ and $\tilde{T}$ commute ($T$ and $\tilde{T}$ are
independent)  and it is not a priori fixed; in fact, $L_q$ is not
uniquely defined (see \cite{FTUV94-21} and \cite{WOZA}).

Consider two $q$-spinors $X$ and $Z$, and their hermitian conjugates
$X^{\dagger}$ and $Z^{\dagger}$, transforming under the coaction of $SL_q(2)$
with the matrices $T$ and $T^{\dagger}$ respectively by
\begin{equation}\label{2a}
\begin{array}{ll}
X'=T X \quad \quad ,& \quad \quad X^{\dagger}\,' = X^{\dagger} T^{\dagger}
\quad ,\\
Z'=T Z  \quad \quad ,& \quad \quad Z^{\dagger}\,' = Z^{\dagger} T^{\dagger}
\quad .
\end{array}
\end{equation}
\noindent
The form of the eqs. (\ref{2b}) is the result of the equations which express
the
commutation relations among the components of the vectors $X,Z$. Since in
principle  $T$ and $T^{\dagger}$ do not commute, we have to allow for
possibly
non-trivial commutation relations among the components of $X$ and
$Z^{\dagger}$.
Thus, the set of commutation relations left invariant  is given by

\begin{equation}\label{2c}
\begin{array}{ll}
R^{(1)} X_{1} X_{2} = q X_{2} X_{1}\quad ,& \quad
Z_{1}^{\dagger} R^{(2)} X_{2} =  X_{2} Z_{1}^{\dagger}\quad ,\\
q Z_{1}^{\dagger} Z_{2}^{\dagger} = Z_{2}^{\dagger} Z_{1}^{\dagger}
R^{(4)}\quad ,& \quad
Z_{2}^{\dagger} R^{(3)} X_{1}  = X_{1} Z_{2}^{\dagger}\quad ,
\end{array}
\end{equation}
\noindent
where $q$ is the deformation parameter of $R^{(1)}$.
The invariance of the first and third equations is
proven as in Sec.2 and the others similarly. In  the second
equation, for example, we can  check that

$$
Z_{1}^{'\dagger} R^{(2)} X_{2}' =
(Z_{1}^{\dagger} T_{1}^{\dagger}) R^{(2)} (T_{2} X_{2})= Z_{1}^{\dagger}
T_{2} R^{(2)} T_{1}^{\dagger} X_{2}
$$
\begin{equation}\label{2d}
=T_{2} Z_{1}^{\dagger} R^{(2)} X_{2} T_{1}^{\dagger} = T_{2} X_{2}
Z_{1}^{\dagger} T_{1}^{\dagger} = X_{2}' Z_{1}^{'\dagger}
\end{equation}

\noindent
using the second equations in (\ref{2b}) and (\ref{2c}), respectively, in the
second and fourth equalities. In particular, if $R^{(2)} = I=
R^{(3)}$,  $T$ and $T^{\dagger}$  commute, which is reflected in
the fact that the components of $X$ and $Z^{\dagger}$ commute.

Let us use the above construction to introduce another  covariant object
which generalizes (with some restrictions)   the concept of twistor
to the $q$-deformed case.
Let $X$ and $Z^{\dagger}$  be $q$-two-vectors ($q$-spinors)
of $SL_{q} (2)$. Tensoring them we introduce the object
\begin{equation}\label{2e}
K \equiv X Z^{\dagger} \qquad (K_{ij}=X_iZ^{\dagger}_j)\quad.
\end{equation}

\noindent
Then, the transformation of $K$ induced by (\ref{2a}) (the $q$-Lorentz
coaction)
is
\begin{equation}\label{2f}
\varphi:\;K \longmapsto K'= TK T^{\dagger} \qquad
(K'_{ij}=T_{im}K_{mn}T^{\dagger}_{nj})\quad.
\end{equation}

\noindent
The entries of $K$ are, of course, non-commuting. We shall see that these
commutation relations can be  expressed by a closed and simple equation
which permits to extract the algebra generated by the entries of $K$
without considering its  explicit realization in terms of the components of $X$
and
$Z^{\dagger}$.
Using the above relations we may now derive the equation describing the
commutation relations which define the algebra generated by the entries of
$K$. With $K_{1} = X_{1} Z_{1}^{\dagger}$
$(K_{1\;ij,kl}=(K \otimes {\bf 1} )_{ij,kl}=X_iZ^{\dagger}_k \delta_{jl})$
and $K_{2}= X_{2} Z_{2}^{\dagger}$
$(K_{2\;ij,kl}=( {\bf 1} \otimes K)_{ij,kl}= \delta_{ik}X_jZ_l^{\dagger})$,
we find using (\ref{2c}) that

\begin{equation}\label{2fa}
\begin{array}{l}
R^{(1)} K_{1} R^{(2)} K_{2} = R^{(1)} X_{1} Z_{1}^{\dagger} R^{(2)} X_{2}
Z_{2}^{\dagger}
=R^{(1)} X_{1} X_{2} Z_{1}^{\dagger} Z_{2}^{\dagger} \\
\qquad = X_{2} X_{1}
Z_{2}^{\dagger} Z_{1}^{\dagger} R^{(4)}=
X_{2} Z_{2}^{\dagger} R^{(3)} X_{1} Z_{1}^{\dagger} R^{(4)}\quad .
\end{array}
\end{equation}

\noindent
Hence, the commuting properties of the quantum twistor are given by

\begin{equation}\label{2g}
R^{(1)} K_{1} R^{(2)} K_{2} = K_{2} R^{(3)} K_{1} R^{(4)} \quad,
\end{equation}

\noindent
which   is  nothing else (see \cite{AKR,PRAGA}) than the
reflection equation (RE) with no spectral parameter
dependence (see \cite{K-SK,K-SA} and references
therein\footnote{Equations of this type were independently used in
\cite{MAJID} in the context of braided algebras.}). Eq. (\ref{2g})
reflects, with (\ref{2b}), the invariance of the commuting properties of the
entries of $K$ by the coaction (\ref{2f}). As shown here, eq. (\ref{2g})
also follows from  interpreting $K$ as an object made out of two $q$-`vectors'
with commuting properties defined by (\ref{2c}).

Let $X$, $Z$ be two $q$-two-vectors (spinors).
We may construct  the following hermitian objects ({\it quantum twistors})
\begin{equation}\label{3d}
K=XX^{\dagger} \quad \mbox{or} \quad K=XZ^{\dagger}+ZX^{\dagger} \quad ,
\end{equation}
\noindent
and find that (\ref{2f})
preserves the hermiticity property of $K$. As we shall see,
the  quantum determinant  ($det_qK$)  of
$K=XX^{\dagger}$  (null $q$-twistor) is necessarily zero
(as it would be as well for $XZ^{\dagger}$). In contrast,
the $q$-twistor  $K=XZ^{\dagger}+ZX^{\dagger}$ has   $det_qK \neq 0$.

To compute the commutation properties of these hermitian $K$ matrices,
the complete set of relations
among  $X$, $Z$, $X^{\dagger}$ and $Z^{\dagger}$ are required. Thus, besides
(\ref{2c}), we require the following set of covariant relations
\begin{equation}\label{tw2}
\begin{array}{ll}
R^{(1)} X_{1} Z_{2} =  Z_{2} X_{1}\quad , & \quad
Z_{1}^{\dagger} R^{(2)} Z_{2} =  Z_{2} Z_{1}^{\dagger}\quad ,\\
X_{1}^{\dagger} Z_{2}^{\dagger} = Z_{2}^{\dagger} X_{1}^{\dagger}
R^{(4)}\quad , & \quad
X_{2}^{\dagger} R^{(3)} X_{1}  = X_{1} X_{2}^{\dagger}\quad ,
\end{array}
\end{equation}
\noindent
the structure of which is again dictated from (\ref{2b}) by covariance.
It is easily seen, using eqs. (\ref{2c}), (\ref{tw2})
and Hecke's condition for the matrix
$R^{(1)}=R^{(4)\,t}$, that the commutation properties
of the $q$-twistors $K$ are also governed by eq. (\ref{2g}).

\vspace{0.3cm}

Notice that $K=XZ^{\dagger}+ZX^{\dagger}$ in (\ref{3d}) is
constructed from two parts, each one
of them satisfying the
same algebra relations (\ref{2g}):
\begin{equation}\label{br1}
\quad K=K^{(1)}+K^{(2)}\quad,\quad
K^{(1)} \equiv XZ^{\dagger} \quad , \quad K^{(2)} \equiv ZX^{\dagger} \quad.
\end{equation}
\noindent
These two pieces have  specific commutation properties among themselves.
Indeed,
the (mixed) commutation relations (\ref{tw2}) lead to the following
non-commuting property
between the matrices $K^{(1)}$ and $K^{(2)}$
(non-symmetric under the interchange of $K^{(1)}$ and $K^{(2)}$)
\begin{equation}\label{br2}
R^{(1)}K^{(1)}_1R^{(2)}K^{(2)}_2=
K^{(2)}_2R^{(3)}K^{(1)}_1({\cal P}R^{(4)}{\cal P})^{-1} \quad.
\end{equation}
\noindent
Setting $R^{(1)}= R_{12}$ or $R_{21}^{-1}$
produces two   different equations
which transform into each other by the  exchange
$K^{(1)}_i \leftrightarrow K^{(2)}_i$.
Both equations had to be possible since $K=K^{(1)}+K^{(2)}$
is symmetric under this exchange.
Eq. (\ref{br2}), which allows the sum (\ref{br1}) of two objects to satisfy
the same commutation properties is an example of `additive
braiding\footnote{For
a discussion of {\it braided geometry} and the role of $q$ see \cite{MB}.}
equation' \cite{MeyerM} here obtained  from the commutation relations
(\ref{tw2}) (`multiplicative braidings' are also possible,
see \cite{MAJID,K-SA}).   Within this terminology,  the `mixed' eqs.
(\ref{tw2}) are the {\it braiding relations for $q$-spinors}.

\vspace{0.3cm}

Supposing that $\hat{R}^{(1)}$=${\cal P}R^{(1)}$ has a spectral
decomposition like (\ref{af3})
with a rank three projector $P_+$ and a rank one projector $P_-$, and
that $det_qT$ (\ref{adet}) and $det_q T^{\dagger}$ are central,
the $q$-determinant of the 2$\times$2 matrix $K$ is given by the expression
\begin{equation}\label{gdet}
(det_qK)P_-= (-q) P_- K_1 \hat{R}^{(3)}K_1 P_- \quad.
\end{equation}
\noindent
When $(det_qT)(det_qT^{\dagger})=1$, $det_qK$ is invariant under the coaction
(\ref{2f}). Using the last eq. in (\ref{2b}) and (\ref{adet})
\begin{equation}\label{det3}
\begin{array}{l}
det_q(TKT^{\dagger})= P_- (T_1K_1T_1^{\dagger})
\hat{R}^{(3)}(T_1K_1T_1^{\dagger})P_- \\
= P_-T_1T_2K_1\hat{R}^{(3)} K_1T_2^{\dagger}T_1^{\dagger}P_- =
(det_qT) \,(det_qK)\,(det_qT^{\dagger}) \;.
\end{array}
\end{equation}
\noindent
Thus, since  $(det_qT)=(det_qT^{\dagger})=1$ we obtain
that $det_q(TKT^{\dagger})=det_qK$.
The centrality of $det_qK$ requires some YBE-like conditions on the
$R^{(i)}$ ($i=1,2,3,4$) matrices in (\ref{2g}).

Using the definition (\ref{gdet}) and the $R$-matrix property
$\hat{R}^{(3)}_{ab,cd}=R^{(3)}_{ba,cd}$, we can compute explicitly the
$q$-determinant of $K$ in the following realizations

\vspace{1\baselineskip}
\noindent
{\bf 1.} For the matrix $K=XZ^{\dagger}$
(and hence for the $q$-twistor $K=XX^{\dagger}$)
\begin{equation}\label{det5}
\begin{array}{l}
(det_qK)P_{-\; ij,kl}= P_{-\; ij,ab}K_{ac}\hat{R}^{(3)}_{cb,mn}
K_{mp}P_{-\;pn,kl}
\propto  \epsilon^q_{ij} \epsilon^q_{ab} X_a Z^{\dagger}_c R^{(3)}_{bc,mn}
X_mZ^{\dagger}_p \epsilon^q_{pn}\epsilon^q_{kl} \\
\hspace{2cm}  =\epsilon^q_{ij} \epsilon^q_{ab} X_a X_b Z^{\dagger}_n
Z^{\dagger}_p \epsilon^q_{pn}\epsilon^q_{kl} =
\epsilon^q_{ij}(X^t\epsilon^qX)(Z^t\epsilon^qZ)^{\dagger}\epsilon^q_{kl} =0
\end{array}
\end{equation}
\noindent
since $(X^t\epsilon^qX)=0=(Z^t\epsilon^qZ)$. This reflects the well-known
fact in non deformed twistor theory that twistors constructed out of
two spinors determine  null length vectors.

\vspace{1\baselineskip}
\noindent
{\bf 2.} For the $q$-twistor $K=XZ^{\dagger}+ ZX^{\dagger}$, a similar calculus
to the previous one gives
\begin{equation}\label{det6}
\begin{array}{ccl}
(det_qK)P_{-\; ij,kl}& \propto &
\epsilon^q_{ij} \epsilon^q_{ab} (X_a Z^{\dagger}_c
+Z_a X^{\dagger}_c) R^{(3)}_{bc,mn}
(X_mZ^{\dagger}_p+Z_m X^{\dagger}_p) \epsilon^q_{pn}\epsilon^q_{kl} \\
\, &=& \epsilon^q_{ij}[(X^t\epsilon^q Z)(X^t\epsilon^q Z )^{\dagger}
+(Z^t\epsilon^q X)(Z^t\epsilon^q X )^{\dagger}]\epsilon^q_{kl} \;
\neq 0 \quad.
\end{array}
\end{equation}
\noindent
Thus, to get twistors with non-null $q$-determinant we need four spinors
in the definition of $K$ (notice that $X$, $Z$, $X^{\dagger}$ and
$Z^{\dagger}$ are all algebraically independent). If  the scalar products
$(X^t\epsilon^q Z)$ and $(Z^t\epsilon^q X )$ are  central
in the algebra generated by  $X$, $Z$, $X^{\dagger}$ and $Z^{\dagger}$
the $q$-determinant of $K$ is also central.

\subsection{An example of $q$-Minkowski space}

$\quad$ Since, by assumption, $T$ and $\tilde{T}=(T^{\dagger})^{-1}$
are $SL_q(2)$ matrices, {\it i.e.},
\begin{equation}\label{3f}
R_{12} T_{1} T_{2} = T_{2} T_{1} R_{12} \quad , \quad
R_{21}T^{\dagger}_{1} T^{\dagger}_{2} = T^{\dagger}_{2} T^{\dagger}_{1}
R_{21}\quad ,
\end{equation}
\noindent
the first and third equations of (\ref{2b})  are fulfilled if
$R^{(1)} = R_{12}$ and $R^{(4)} = R_{21}$.
Assuming, for instance,  the specific non-trivial commutation
relations between $T$ and $T^{\dagger}$ given by
\begin{equation}\label{3w}
T^{\dagger}_2R_{12} T_{1} = T_{1} R_{12}T^{\dagger}_2
\end{equation}
\noindent
we find that $R^{(2)} = R_{21}$. Then, eq.
(\ref{2g}) leads to \cite{AKR}
\begin{equation}\label{3x}
R_{12} K_{1} R_{21} K_{2} = K_{2} R_{12} K_{1} R_{21}
\end{equation}
\noindent
In terms of the entries of $K = \left( \begin{array}{cc}
                              \alpha & \beta \\
                               \gamma & \delta
                             \end{array} \right)$
eq. (\ref{3x}) leads to the algebra (see also \cite{MAJID,K-SA})
\begin{equation}\label{algebra}
\begin{array}{lll}
\alpha \beta = q^{-2} \beta \alpha \quad, &
\alpha \gamma = q^{2} \gamma \alpha \quad, &
[ \beta, \gamma ] = (\lambda/q)( \alpha \delta - \alpha^2) \quad, \\
\alpha \delta =  \delta \alpha  \quad, &
 [ \beta, \delta ] = -(\lambda/q) \alpha \beta \quad, &
 [ \gamma , \delta] = (\lambda/q) \gamma \alpha \quad.
\end{array}
\end{equation}
\noindent
These commutation properties are preserved by (\ref{2f}) and define the
{\it quantum
Minkowski algebra} ${\cal M}_q$  of \cite{WSSW,SWZ,OSWZ}. Its
linear central term (the $q$-trace of $K$ [(\ref{qtr})])
is identified with time,
\begin{equation}\label{xo}
x^0 \sim tr_qK=q^{-1} \alpha + q \delta \quad,
\end{equation}
\noindent
and the $q$-determinant, defined by (\ref{gdet}) with $\hat{R}^{(3)}
=\hat{R}$,
\begin{equation}\label{3z}
(det_qK ) P_-= (-q) P_-K_1 \hat{R} K_1P_-
= (\alpha \delta -q^2 \beta \gamma ) P_- \quad,
\end{equation}
\noindent
gives the invariant quadratic central element which is identified
with the  $q$-{\it Minkowski length}. For other examples of $q$-Minkowski
algebras and further
discussions, see \cite{FTUV94-21}.

\setcounter{equation}{0}

\section{\bf Non-commutative differential calculus}

$\quad$ The development of a non-commutative differential calculus
on the quantum groups (see, {\it e.g.}, \cite{WOR1,JURCO,ZU-MPX}
and references therein and \cite{QGR} for reviews) as well as on the
quantum spaces  (see, {\it e.g.}, \cite{WZ,MANIN2,Pusz-WO}
and references therein)
requires including, in a first stage,  derivatives and differentials.
A Cartan calculus involving the inner derivation and the Lie
derivative (see \cite{CARTAN,WATTS} and references therein)
and the $q$-Hodge star operator \cite{MEQ} may also be introduced
but will not be  discussed here.

\subsection{ Differential calculus {\it \`a la } Wess-Zumino}

$\quad$ The covariant differential calculus on the quantum planes was discussed
by J. Wess and B. Zumino \cite{WZ} for the $A_l$-type plane in which the
corresponding $\hat{R}$-matrix has two different eigenvalues. The method
developed there was generalized to the quantum orthogonal plane
($B_l$, $D_l$-type)
\cite{WA-ZPC49,SONG2,OGI}, to the $q$-symplectic quantum plane ($C_l$-type)
\cite{SONG2}, and, in general, to a  quantum plane given in terms of an
$\hat{R}$-matrix with any number of eigenvalues \cite{ZU-SA}. This approach
is also useful to discuss the differential calculus on the $q$-Minkowski space
\cite{OSWZ,SONG}.

\subsubsection{The formalism}

$\quad$ Consider the associative algebra ({\it quantum space}) generated by
variables
(`coordinates') $x^i$, $i=1,2,...,n$, with relations (cf. (\ref{1d}))
\begin{equation}\label{4a}
x^ix^j-B_{ij,kl}x^kx^l=0
\end{equation}
\noindent
which may be rewritten  $X=(x^1,..., x^n)$
\begin{equation}\label{4b}
(I-B_{12})X_1X_2=0 \quad.
\end{equation}
\noindent
Following \cite{WZ} (see also
\cite{MANIN2,ZU-MPX}), we first introduce the exterior derivative $d$.
$d\, : \, x^i \mapsto dx^i$ is  nilpotent,  $d^2=0$,
and satisfies  Leibniz's  rule (for a modified version see \cite{FP}).
In terms of the derivatives,
\begin{equation}\label{4f}
d= dx^i \partial_i \;, \quad  i=1,2,...,n\;, \quad \quad \partial_i \equiv
\frac{ \partial}{ \partial x^i} \quad.
\end{equation}
\noindent
The differential calculus on the quantum space is defined by a set of
quadratic algebraic relations among all the fundamental objects: coordinates
$x^i$,
$q$-differentials $dx^i$ and $q$-derivatives $\partial_i$. These
commutation relations must be introduced according to two essential
requirements:

\noindent
{\bf 1.} Covariance under the transformation (coaction) generated by a
quantum group  matrix  $T$
\begin{equation}\label{4g}
X \mapsto X'=TX \quad , \quad dX \mapsto dX'=T dX \quad , \quad
\partial \mapsto \partial'= (T^t)^{-1} \partial  \quad,
\end{equation}
\noindent
in this way, $d$ is invariant under this   transformation, $d=dX^t \partial=
dX'^t \partial'$.

\noindent
{\bf 2.} Consistency, which means that the quadratic algebra generated by
$x^i$, $dx^i$ and $\partial_i$ ($i=1,2,...,n$) is associative and
there are not independent higher order relations.

\vspace{0.3cm}

In this way, the cross commutation relations among the differential quantities
are given in matrix form  by \cite{WZ}
\begin{equation}\label{4h}
\begin{array}{ll}
(I-B_{12})X_1X_2 =0 \quad, & \;  \\
\partial_2 \partial_1 (I- F_{12})=0 \quad , & \qquad \partial_1 X_1 =
I + X_2 C_{12}^t  \partial_2 \;, \\
(I+C_{12})dX_1dX_2=0 \quad, & \qquad \partial_1 dX_1
= dX_2 D_{12}^t \partial_2 \;, \\
X_1dX_2=C_{12} dX_1 X_2 \quad , & \;
\end{array}
\end{equation}
\noindent
where $B$, $C$, $D$  and $F$ are numerical matrices to be determined.
The second requirement leads to a set consistency conditions which can be
summarized in the following relations
\begin{equation}\label{4i}
\begin{array}{l}
(I-B_{12})(I+C_{12})=0 \quad, \qquad (I+C_{12})(I-F_{12})=0 \quad,\\
(I_{12} - B_{12})C_{23}C_{12}=C_{23}C_{12}(I_{23}-B_{23}) \quad, \\
C_{12}C_{23}(I_{12}-F_{12}) =(I_{23}-F_{23})C_{12}C_{23}\quad, \\
C_{12}C_{23}D_{12}=D_{23}C_{12}C_{23}\quad, \qquad D=C^{-1} \quad.
\end{array}
\end{equation}
\noindent
Now, for a given $\hat{R}$-matrix, the solution of these consistency conditions
may be obtained by writing, in a appropriate way,
the matrices $(I-B)$, $(I-F)$  and $(I+C)$ in terms of the projectors
obtained from the spectral decomposition of $\hat{R}$.

\subsubsection{The two-dimensional quantum plane case}

$\quad$ The $GL_q(2)$ $\hat{R}$-matrix has the spectral
decomposition (\ref{af3}) where the projector
$P_+$ ($P_-$) is the deformed version of the symmetrizer (antisymmetrizer)
in two dimensions. In view of (\ref{proj}), it is not surprising that
these  projectors also allow  us  to define quantum planes.
For instance, as for $q$=1,
the commutation relations for the two-dimensional quantum space are obtained by
requiring the vanishing of their $q$-antisymmetric products
\begin{equation}\label{4j}
P_-X_1X_2=0 \quad.
\end{equation}
\noindent
This equation is obtained from eq. (\ref{1d})   since,
using there  the spectral decomposition (\ref{af3})
of $\hat{R}$  one gets
$(qP_+-q^{-1}P_-)X_1X_2=q(P_++P_-)X_1X_2$.
Analogously, substituting $dX$ for $\Omega$ in (\ref{1d1})  it follows that
the $q$-symmetric  products of $q$-differentials  must vanish,
\begin{equation}\label{4k}
P_+dX_1dX_2=0 \quad
\end{equation}
\noindent
since, using (\ref{af3}), eq.   (\ref{1d1}) may be rewritten as
$(qP_+-q^{-1}P_-) dX_1 dX_2 = -q^{-1} (P_++P_-)dX_1dX_2$.
Obviously, all these relations involving  projectors are preserved by the
$GL_q(2)$ coaction.
The consistency equations (\ref{4i}) are fulfilled by taking \cite{WZ}
\begin{equation}\label{4l}
\begin{array}{rl}
I-B=I-F=P_- & \;  \longrightarrow  \quad B=F= q^{-1} \hat{R} \quad, \\
I+C=P_+ & \; \longrightarrow  \quad C=q \hat{R} \quad,
\end{array}
\end{equation}
\noindent
and the covariant differential calculus for $q$-two-vectors
(or for $q$-spinors) is defined by (compare with the matrix notation used
in (\ref{4h}))
\begin{equation}\label{4m}
\begin{array}{ll}
x^i dx^j = q \hat{R}_{ij,kl}dx^kx^l \quad , &
dx^i dx^j = -q \hat{R}_{ij,kl}dx^kdx^l \quad , \\
\partial_k x^i = \delta^i_k + q \hat{R}_{ij,kl}x^l \partial_j  \quad , &
\partial_i \partial_j =q^{-1} \hat{R}_{lk,ji}\partial_k \partial_l \quad, \\
\partial_k dx^i = q^{-1}  \hat{R}_{ij,kl}^{-1} dx^l \partial_j  \quad . & \,
\end{array}
\end{equation}
\noindent
In the limit $q$=1, $\hat{R}_{ij,kl}= \delta_{il} \delta_{jk}$  reproduces
the usual relations.
All these relations are preserved under the quantum group transformations
(\ref{4g}).
For instance, multiplying the inhomogeneous eq. $\partial_1 X_1 =
I + q X_2 \hat{R}_{12}  \partial_2$
by $(T_1^t)^{-1}$ from the left and
by $T_1^{t}$ from the right and using the RTT relation in the form
$(T_1^t)^{-1} \hat{R}_{12}T_1^t=T_2^t \hat{R}_{12}(T_2^t)^{-1}$
(cf. (\ref{hat})), it follows that  $\partial_1' X_1' =
I +q  X_2' \hat{R}_{12}  \partial_2'$.

\subsubsection{Differential calculus on $q$-Minkowski space }

$\quad$ The previous formalism was used in \cite{WA-ZPC49}
(see also \cite{SONG2,OGI})
to develop the $SO_q(n)$-covariant differential calculus on $n$-dimensional
$q$-Euclidean spaces. In this case,  the  $\hat{R}$-matrix  of $SO_q(n)$,
has three eigenvalues \cite{FRT}. In
contrast, the situation for the $q$-Lorentz group and the $q$-Minkowski
space ${\cal M}_q$ is slightly more complicated.  The
commutation relations among the components of the $q$-Minkowski vector
are computed by using the commutation relations among $q$-spinors
(see \cite{WSSW,OSWZ,SONG} for explicit calculations and formulas). This leads
to two different $\hat{R}$-matrices, $\hat{{\cal R}}_I$ and
$\hat{{\cal R}}_{II}$, both 16$\times$16 matrices satisfying  the Yang-Baxter
equation and with eigenvalues $q^2$, $q^{-2}$  and  $-$1;
$\hat{{\cal R}}_I$ and  $\hat{{\cal R}}_{II}$ give rise to three projectors
each. In each case, one of the
subspaces can be further decomposed so that there are four independent
projectors in all. These projectors allow us to write down
commutation  relations for the $q$-symmetric variables (coordinates)
and for the $q$-antisymmetric ones ($q$-differentials). Thus, as it
was shown in \cite{OSWZ}, the commutation relations among coordinates,
differentials and
derivatives satisfying the requirements of covariance and consistency, which
define the differential calculus on the $q$-Minkowski space, can be
expressed in terms of $\hat{{\cal R}}_I$ and  $\hat{{\cal R}}_{II}$.

\subsection{ RE formalism and $q$-Minkowski space calculus}

$\quad$ We shall now look at the $q$-Minkowski space  differential
calculus  by expressing the different commutation relations
in terms of  appropriate RE, so that the $q$-derivatives (${\cal D}_q$) and
the $q$-forms ($\Lambda_q$)  algebras will also be defined by RE
\cite{AKR,FTUV94-21}.
Consider first an object $Y$ transforming covariantly {\it i.e.},
\begin{equation}\label{ta}
Y \longmapsto Y' = (T^{\dagger})^{-1} Y T^{-1} \quad
\end{equation}
\noindent
(cf. (\ref{2f}), which will be taken as contravariant).
It is easily seen, using (\ref{2b}), that
the invariance of the commutation properties of the matrix elements of $Y$
 gives
\begin{equation}\label{yalg}
R^{(1)} Y_{1} R^{(3)\,-1} Y_{2} = Y_{2} R^{(2)\,-1} Y_{1} R^{(4)}
\end{equation}
\noindent
A quadratic  and $L_q$-invariant element ($q$-determinant)
is defined through
\begin{equation}\label{tf}
(det_{q} Y) P_{-} =(-q^{-1}) P_{-} Y_{1} \hat{R}^{(3)\,-1} Y_{1} P_{-} \quad;
\end{equation}
\noindent
${\Box}_q \equiv det_qY$ becomes  the $q$-{\it D'Alembertian}
once the components of $Y$ are associated with the $q$-derivatives.
As the $K$ matrix entries were associated with the generators of
${\cal M}_q$, we shall consider the elements of $Y$ as generating
the algebra ${\cal D}_q$
of the $q$-Minkowski derivatives. For the example of Sec.{\bf 3.2}
the commutation properties are given by
\begin{equation}\label{yalg2}
R_{12} Y_{1} R_{12}^{-1} Y_{2} = Y_{2} R_{21}^{-1} Y_{1} R_{21} \quad, \quad
Y = \left( \begin{array}{cc}
                              u & v \\
                              w & z
                             \end{array} \right) \quad,
\end{equation}
\noindent
and the $q$-determinant
\begin{equation}\label{tf11}
(det_{q} Y) P_{-} =(-q^{-1}) P_{-} Y_{1} \hat{R}^{-1} Y_{1} P_{-}
= (uz - q^{-2} v w) P_-\quad
\end{equation}
\noindent
is central, $[ \Box_q , Y]=0$.
We now need to establish the commutation
properties among coordinates
and derivatives extending   the classical
relation $\partial_{\mu} x^{\nu} = x^{\nu} \partial_{\mu} + \delta_{\mu}^{\nu}
$, $\partial^{\dagger}$=$- \partial$ to the non-commutative  case,
in a $L_q$-invariant manner.
This is achieved  by
 an {\it inhomogeneous} RE \cite{AKR} of the form
\begin{equation}\label{tk}
Y_{2} R^{(1)} K_{1} R^{(2)} = R^{(3)} K_{1} R^{(1)\,-1} Y_{2} + \eta J,
\end{equation}

\noindent
where $\eta$ is a constant, $\eta J$$\rightarrow$$I_{4}$ in the
$q$$ \rightarrow$$1$ limit, and $J$ is invariant,
\begin{equation}\label{tl}
J \longmapsto (T_{2}^{\dagger})^{-1} T_{1} J T_{1}^{\dagger} T_{2}^{-1} = J
\,, \qquad  T_{1} J T_{1}^{\dagger}  =  T_{2}^{\dagger}J T_{2} \,.
\end{equation}
\noindent
As for $J$,  setting $J\equiv  J' {\cal P}$ in eq. (\ref{tl}) gives
$T_{1} J' T_{2}^{\dagger}  =  T_{2}^{\dagger}J' T_{1}$,
hence $J = R^{(3)} {\cal P}$ (the same
result follows if  we set $J= {\cal P} J'$). In the previous example
$J = R^{(3)} {\cal P} = R_{12} {\cal P}$, and
in order to have the inhomogeneous term in the simplest form (the analogue of
the $\delta_{\nu}^{\mu}$ of the $q=1$ case) it is convenient to take $\eta
= q^{2}$. Then, the
commutation relations of the entries of $K$ (generators of the algebra
${(\cal M}_q$ of coordinates) and those of $Y$ (generators of the algebra
${(\cal D}_q$ of derivatives) are given by
\begin{equation}\label{ts}
Y_{2} R_{12} K_{1} R_{21} = R_{12} K_{1} R^{-1}_{12} Y_{2} + q^2
R_{12}{\cal P} \quad.
\end{equation}
\noindent
This equation is not invariant under hermitian conjugation. In fact, it is not
possible to have {\it simultaneously} coordinates and derivatives with the
usual hermiticity properties ($K=K^{\dagger}$ and $Y=-Y^{\dagger}$)
\cite{OSWZ} (see also \cite{OZ}).

The determination of the commutation relations for the $q$-De Rham complex
requires incorporating the exterior derivative $d$. To the  four generators
of ${\cal M}_{q}$ and of ${\cal D}_{q}$ we now add the four elements  of $dK$
($q$-one-forms), which generate the
de Rham complex algebra  $\Lambda_{q}$
(the degree of a form is defined as in the classical case).
As in (\ref{4g}), $d$ commutes with the  $q$-Lorentz coaction (\ref{2f}), so
that
\begin{equation}\label{tt}
dK' = T d K T^{\dagger} \quad .
\end{equation}

\noindent
Applying $d$ to  (\ref{3x}) we obtain
\begin{equation}\label{tu}
R_{12} dK_{1} R_{21} K_{2} + R_{12} K_{1} R_{21} d K_{2} = d K_{2} R_{12}
K_{1}R_{21} + K_{2} R_{12} d K_{1} R_{21} \quad .
\end{equation}

\noindent
We now use that $R_{12} = R^{-1}_{21} + \lambda {\cal P}$
(and the same for 1$\leftrightarrow$2) to replace one $R$
in each term  in such a way that the terms in ${\cal P}K_{1} R_{21} dK_{2}$
and in ${\cal P}dK_1R_{21}K_2$ may
cancel. In this way we obtain two solutions of (\ref{tu}). Since the
relations obtained are not  invariant under
hermitian conjugation, we may use one of them
for $dK$ and the other for the hermitian conjugate
$dK^{\dagger} \equiv (dK)^{\dagger}$
\begin{equation}\label{tv}
R_{12} K_{1} R_{21} d K_{2} = d K_{2} R_{12} K_{1} R^{-1}_{12} \quad,\quad
R_{12} d K^{\dagger}_{1} R_{21} K_{2} = K_{2} R_{12} dK^{\dagger}_{1}
R_{12}^{-1}  \;,
\end{equation}

\noindent
from which follows that
\begin{equation}\label{tw}
R_{12}dK_{1} R_{21} d K_{2} = -d K_{2} R_{12} dK_{1} R^{-1}_{12} \quad,\quad
R_{12} d K^{\dagger}_{1} R_{21} dK_{2} = -dK_{2}
R_{12} dK^{\dagger}_{1}  R_{12}^{-1}  \;.
\end{equation}
\noindent
 The exterior derivative $d$,
as its invariance suggests,  has  the form
\begin{equation}\label{tac}
d=  tr_{q} (dKY)\;.
\end{equation}
\noindent
For an explicit comparison with the formalism of Sec.{\bf 4.1.3} \cite{OSWZ}
see
\cite{FTUV94-21,KARPACZ}.

What about the physical applications of non-commutative geometry to physics?
One of the reasons for introducing $q$ was to see whether the infinities in
quantum field theory could be made milder. For the moment, however, there is
no `$q$-special relativity theory' or `$q$-deformed quantum field theory'.
We shall conclude by mentioning just one interesting
application to particle physics,
the Connes-Lott version of the standard model \cite{CL}; very recently,
it has been used to give an indication of the Higgs mass \cite{KS}.

\vspace{1\baselineskip}

\noindent
{\bf Acknowledgements}: The authors wish to thank P.P. Kulish for helpful
discussions.
This research has been  partially supported by
a CICYT research grant.

\appendix

\section{Appendix. Notation and useful expressions}

$\quad$ We list here some expressions and conventions that are useful in the
main
text. `$RTT$' relations as those in (\ref{uc}), (\ref{2b}) follow the
usual conventions
i.e., the 4$ \times $4 matrices\footnote{Most  formulae in  this Appendix
are also valid for the general $GL_q(n)$ case by
setting $i,j,k,...=1,2,...,n$.} $T_1$, $T_2$ are the tensor products
\begin{equation}\label{aa}
T_{1} = T \otimes I \quad , \quad T_{2} = I \otimes T \quad .
\end{equation}

\noindent
The tensor product of two matrices, $C= A \otimes B$, reads in components
\begin{equation}\label{ab}
C_{ij,kl} = A_{ik} B_{jl} \quad ,
\end{equation}

\noindent
so that the comma separates the row and column indices of the two matrices.
Thus, $(A_{1})_{ij,kl} = A_{ik} \delta_{jl}\,$; $(A_{2})_{ij,kl}
= A_{jl} \delta_{ik}\,$. The
transposition in the first and second spaces is given by
\begin{equation}\label{ac}
C^{t_{1}}_{ij,kl} = C_{kj,il} \quad , \quad C^{t_{2}}_{ij,kl} = C_{il,kj}
\quad, \end{equation}

\noindent
i.e., $C^{t_{1}} = A^{t} \otimes B$ (resp. $C^{t_{2}}= A \otimes
B^{t}$) is given by a matrix in which the blocks 12 and 21 are
interchanged (each of the four blocks is replaced by its transpose). Of course,
$C_{ij,kl}^{t_{1}t_{2}} \equiv  C^{t}_{ij,kl} = C_{kl,ij}$ is the ordinary
transposition. Similarly, the traces in the first and second spaces are given
by
\begin{equation}\label{ad}
(tr_{(1)} C)_{jl} = C_{ij,il} \quad , \quad (tr_{(2)} C)_{ik} = C_{ij,kj}
\quad .
\end{equation}

\noindent
They correspond, respectively, to replacing the $4 \times 4$ matrix
$C$ by the $2 \times 2$ matrix resulting from
adding its two diagonal boxes or by the $2 \times 2$
matrix obtained by taking the trace of
each of its four boxes. If $\,C= A \otimes B$, $tr_{(1)}C= (tr A)B$ and
$tr_{(2)} C= A(tr B)$.

The action of the permutation matrix ${\cal P}_{12} \equiv {\cal P} $
is defined by $({\cal P} C {\cal P})_{ij,kl}=C_{ji,lk}$
(${\cal P} (A \otimes B) {\cal P} = B \otimes A$ if the entries
of $A$ and $B$ commute); thus
\begin{equation}\label{ae}
({\cal P} A_{1} {\cal P})_{ij,kl} = (A_1)_{ji,lk}= A_{jl}
\delta_{ik} = (A_{2})_{ij,kl} \quad ;
\end{equation}

\noindent
$({\cal P} C)_{ij,kl} = C_{ji,kl} \;,\; (C {\cal P})_{ij,kl} =
C_{ij,lk}$. Explicitly, ${\cal P}$=${\cal P}^{-1}$ is given by
\begin{equation}\label{af}
{\cal P}=
\left[
\begin{array}{llll}
1 & \, & \, & \, \\
\, & 0 & 1 & \, \\
\, & 1 & 0 & \, \\
\, & \, & \, & 1
\end{array}
\right]
\quad , \quad {\cal P}_{ij,kl}= \delta_{il} \delta_{jk} \quad;
\end{equation}

\noindent
acting from the left (right) it interchanges the second and third rows
(columns).

   For $GL_q(2)$ (and $SL_q(2)$), the $R_{12} (q) \equiv R_{12}
\equiv R$ and ${\cal P} R_{12} \equiv
\hat{R}_{12}\equiv\hat{R}$ matrices are given by
\begin{equation}\label{af2}
R=
\left[
\begin{array}{llll}
q & \, & \, & \, \\
\, & 1 & 0 & \, \\
\, & \lambda & 1 & \, \\
\, & \, & \, & q
\end{array}
\right]\;,
\quad
 \hat{R}=\left[
\begin{array}{llll}
q & \, & \, & \, \\
\, & \lambda & 1 & \, \\
\, & 1 & 0 & \, \\
\, & \, & \, & q
\end{array}
\right] = \hat{R}^t \;,
\end{equation}
\begin{equation}\label{ag}
R_{12} (q^{-1})= R_{12}^{-1} (q) \;,\qquad
\hat{R}_{12}^{-1} (q)= \hat{R}_{21} (q^{-1}) \; ;
\end{equation}

\noindent
where $\lambda \equiv q - q^{-1}$;
$\hat{R}_{21}= {\cal P} \hat{R}_{12} {\cal P}$. In terms of $\hat{R}$,
the RTT equation reads
\begin{equation}\label{hat}
\hat{R}_{12}T_1T_2=T_1T_2\hat{R}_{12} \quad , \quad
(\hat{R}_{ij,ab}T_{ak}T_{bl}=T_{ic}T_{jd}\hat{R}_{cd,kl}) \quad .
\end{equation}
\noindent
Similarly,
${\cal P} R_{12} {\cal P}=R_{21}=R_{12}^t$, but the last equality is due to the
specific form of $R_{12}$.
$\hat{R}$ satisfies Hecke's condition
\begin{equation}\label{ae2}
\hat{R}^{2} - \lambda \hat{R} - I = 0 \quad , \quad (\hat{R} - q)
(\hat{R} + q^{-1}) =0
\end{equation}

\noindent
and
\begin{equation}\label{af3}
\hat{R} = q P_{+} -q^{-1} P_{-} \;, \quad \hat{R}^{-1} = q^{-1} P_{+}
- q P_{-}\;,\quad
[\hat{R} , P_{\pm} ] =0\; , \quad P_{\pm} \hat{R} P_{\mp} = 0 \;,
\end{equation}

\noindent
where  the projectors $P_{\pm \,12} \equiv P_{\pm}$ ($q$-(anti)symmetrizer)
are given by
\begin{equation}\label{ag2}
P_{+} = \frac{1}{[2]}
\left[
\begin{array}{cccc}
[2] & 0 & 0 & 0\\
0 & q & 1 & 0\\
0 & 1 & q^{-1} & 0\\
0 & 0 & 0 & [2]
\end{array}
\right] \;,
\qquad
P_{-}=\frac{1}{[2]}
\left[
\begin{array}{cccc}
0 & 0 & 0 & 0\\
0 & q^{-1} & -1 & 0\\
0 & -1 & q & 0\\
0 & 0 & 0 & 0
\end{array}
\right]\;,
\end{equation}

\noindent
with $[2] \equiv ( q + q^{-1})$. It is often convenient
to express the 4$\times$4 matrix $P_{-}$ in the form
\begin{equation}\label{ah}
(P_{-})_{ij,kl} = \frac{1}{[2]} \epsilon_{ij}^{q} \epsilon^{q}_{kl} \;,\qquad
([x] \equiv  \frac{q^{x} - q^{-x}}{q - q^{-1}})  \quad,
\end{equation}

\noindent
where $\epsilon^{q}$$=$$-(\epsilon^{q})^{-1} $$\neq $$
(\epsilon^{q})^{t}$  is given in (\ref{epsilon1}).
 The determinant of an ordinary $2 \times 2$ matrix  may be defined
as the proportionality
coefficient in  $(det T)P_{-} = P_{-} T_{1} T_{2}$ where $P_{-}$ is obtained
from (\ref{ag2}) setting $q$=1. The analogous  definition in  the
$q \neq 1$ case
\begin{equation}\label{adet}
(det_q T)P_{-} := P_{-} T_{1} T_{2}\quad, \qquad
(det_q T^{\dagger})P_{-} = T_{2}^{\dagger} T_{1}^{\dagger}P_-\;,
\end{equation}
\noindent
($det_q T^{\dagger}= (det_q T)^{\dagger}$)
leads to the  expression for $det_q T$ given in (\ref{ub}).
For the $K$ matrix, the definition of $det_qK$ is
given by (\ref{gdet}).

The $q$-trace of a matrix $B$ is defined by
(see \cite{FRT,ZU-MPX})
\begin{equation}\label{gtr}
tr_q(B)=tr(DB) \quad , \quad
D=q^2tr_{(2)}({\cal P}(( (R^{(1)})^{t_1})^{-1})^{t_1}) \quad,
\end{equation}
\noindent
where the superscript $t_1$ means transposition in the first space, the trace
$tr_{(2)}$ is taken in the second space  and  $q$ is the deformation
parameter in $R^{(1)}$ (see (\ref{2b}) and  (\ref{2g})).
This $q$-trace is invariant under the quantum group coaction
$B \mapsto TBT^{-1}$
(as well as under the coaction $C \mapsto (T^{\dagger})^{-1} C T^{\dagger}$
since $R^{(1)}$=$R^{(4)\,t}$=${\cal P}R^{(4)}{\cal P}$).
In particular, if $R^{(1)}=R_{12}$ given in (\ref{af2}) the $q$-trace of $B$
is (\ref{xo}),
\begin{equation}\label{qtr}
tr_qB=tr(DB) \quad, \qquad D= \left( \begin{array}{cc}
                                          q^{-1} & 0 \\
                                            0  &  q
                                      \end{array} \right) \quad.
\end{equation}


\begin{thebibliography}{99}

\bibitem{CONNES} A. Connes, {\it G\'eom\'etrie non commutative},
InterEditions (1990)

\bibitem{COQUE} R. Coquereaux, J. Geom. Phys. {\bf 6}, 425 (1989)

\bibitem{COL} A collection of pioneering papers may be found in
{\it Yang-Baxter equation in integrable systems}, M. Jimbo ed.
(Adv. Series in Math. Physics {\bf 10}), World Sci. (1989)

\bibitem{DRINFELD} V.G. Drinfel'd, in Proc. of the 1986 {\it Int. Congr. of
Mathematicians}, MSRI Berkeley, vol. {\bf I}, 798 (1987)

\bibitem{JIMBO} M. Jimbo, Lett. Math. Phys. {\bf 10}, 63 (1985);
ibid. {\bf 11},  247 (1986)

\bibitem{FRT} L.D. Faddeev, N. Yu. Reshetikhin and L.A. Takhtajan, Alg. i Anal.
{\bf 1}, 178 (1989) (Leningrad Math. J. {\bf 1}, 193 (1990))

\bibitem{MANIN} Yu. I. Manin, Commun. Math. Phys. {\bf 123}, 163 (1989);
{\it Topics in non-commutative geometry}, Princeton Univ. Press. (1991)


\bibitem{WZ} J. Wess and B. Zumino, Nucl. Phys. (Proc. Suppl.) {\bf 18B},
302 (1990)


\bibitem{DAVID} E. Corrigan, D.B. Fairlie, P. Fletcher and R. Sasaki,
J. Math. Phys. {\bf 31}, 776 (1990)


\bibitem{ABE} E. Abe, {\it Hopf algebras}, Cambridge Univ. Press (1977)


\bibitem{VZW}S.P. Vokos, B. Zumino  and J. Wess, Z. Phys. {\bf C48}, 65 (1990)

\bibitem{YB} M. Jimbo, Int. J. Mod. Phys. {\bf A4}, 3759 (1989);
H. de Vega, Int. J. Mod. Phys. {\bf A4}, 2371 (1989)

\bibitem{TAK} L. A. Takhtajan, {\it Lectures in quantum groups}, in Nankai
Lectures
in Math. Phys. (M.-L. Ge and B.-H. Zhao eds.),
World Sci. (1990), p. 69;
A. Pressley and V. Chari, Nucl. Phys. (Proc. Suppl.) {\bf 18A},
207 (1990);
H. D. Doebner, J. D. Henning and W. L\"ucke, {\it
Mathematical guide to quantum groups},  in Proc. of the Clausthal
Int. Workshop on Math. Phys., H.-D.
Doebner and J.-D. Henning Eds., Springer-Verlag (1990), p.29;
S. Majid, {\it Foundations of quantum group theory}, Cambridge
Univ. Press (1995)

\bibitem{MANIN2} Yu. I. Manin, Theor. and Math. Phys. {\bf 92}, 997 (1993)

\bibitem{PRAGA} J.A. de Azc\'{a}rraga, P.P. Kulish and F. Rodenas,
{\it Non-commutative geometry and covariance: from the quantum plane to
quantum tensors}, to appear in the Proc.
of the 3rd Colloquium `Quantum groups and Physics' (Czechoslovak J. Phys.,
1995)

\bibitem{POWO} P. Podles  and S. Woronowicz, Commun. Math. Phys. {\bf 130},
381 (1990)

\bibitem{WSSW} U. Carow-Watamura, M. Schlieker, M. Scholl and S. Watamura, Z.
Phys. {\bf C48}, 159 (1990); Int. J. Mod. Phys. {\bf A6}, 3081 (1991)

\bibitem{SWZ} W. Schmidke, J. Wess and B. Zumino, Z. Phys. {\bf C52}, 471
(1991)

\bibitem{OSWZ} O. Ogievetsky, W. B. Schmidke, J. Wess and B. Zumino, Commun.
Math. Phys. {\bf 150}, 495 (1992)

\bibitem{FTUV94-21} J.A. de Azc\'{a}rraga, P.P. Kulish and F. Rodenas,
{\it Quantum groups and deformed special relativity}, FTUV-94-21/IFIC-94-19
(April 1994, hep-th 9405161)

\bibitem{WOZA} S.L. Woronowicz and S. Zakrzewski, Compositio Math.
{\bf 94}, 211 (1994)

\bibitem{AKR} J.A. de Azc\'{a}rraga, P.P. Kulish and F. Rodenas,
Lett. Math. Phys. {\bf 32}, 173 (1994)

\bibitem{K-SK} P.P. Kulish and E.K. Sklyanin, J. Phys. {\bf A25}, 5963 (1992)

\bibitem{K-SA} P.P. Kulish and R. Sasaki, Progr. Theor. Phys. {\bf 89}, 741
(1993)

\bibitem{MAJID} S.  Majid, in {\it Quantum Groups}, Lect. Notes  Math. {\bf
1510},
79 (1992); J. Math. Phys. {\bf 32}, 3246 (1991); ibid {\bf 34}, 1176 (1993)

\bibitem{MB} S. Majid, {\it Introduction to braided geometry and
$q$-Minkowski space}, DAMTP 94-68 (Varenna lectures, July 1994)

\bibitem{MeyerM}  U. Meyer, {\it The $q$-Lorentz group and braided
coaddition on $q$-Minkowski space}, DAMTP 93-45 (revised, Jan. 1994);
S. Majid, J. Math. Phys. {\bf 34}, 2045 (1993);
S. Majid and U. Meyer, Z. Phys. {\bf C63}, 357 (1994)

\bibitem{WOR1} S. L. Woronowicz,  Publ. RIMS Kyoto Univ. {\bf 23},
117 (1987); Commun. Math. Phys. {\bf 122}, 125 (1989)

\bibitem{JURCO} B. Jurco, Lett. Math. Phys. {\bf 22}, 177 (1991)

\bibitem{ZU-MPX} B. Zumino, {\it Introduction to the differential geometry
of  quantum groups},
in  {\it Mathematical Physics X},  K.  Schm\"udgen ed.,
Springer-Verlag (1992), p. 20

\bibitem{QGR} O. Bernard, Progr. Theor. Phys. Suppl. {\bf 102}, 49 (1990);
F. M\"uller-Hoissen, Nucl. Phys. {\bf 25}, 1703 (1992);
P. Aschieri and L. Castellani, Int. J. Mod. Phys. {\bf A8}, 1667 (1993)

\bibitem{Pusz-WO} W. Pusz and S.L. Woronowicz, Rep. Math. Phys. {\bf 27},
231 (1990)

\bibitem{CARTAN} C. Chryssomalakos, P. Schupp and B. Zumino, {\it
Induced extended calculus on the quantum plane}, LBL-35034, hep-th 9401141

\bibitem{WATTS} P. Watts, {\it Differential geometry on Hopf algebras
and quantum groups}, LBL-36537/UCB-TH 94/35

\bibitem{MEQ} U. Meyer, {\it Wave equations on $q$-Minkowski spaces},
DAMTP 94-10 (1994)

\bibitem{FP} L. D. Faddeev and P. N. Pyatov, {\it The differential calculus
on quantum linear groups},  preprint (1993), hep-th 9402070

\bibitem{WA-ZPC49} U. Watamura, M. Schlieker and S. Watamura,
Z. Phys. {\bf C49}, 439 (1991)

\bibitem{SONG2} X-C. Song and L. Liao, J. Phys. {\bf A25}, 623 (1992)

\bibitem{OGI} O. Ogievetsky, Lett. Math. Phys. {\bf 24}, 245 (1992)

\bibitem{ZU-SA} B. Zumino, Anales de F\'{\i}sica: Monograf\'{\i}as {\bf 1},
vol I, pp. 41 (1993)

\bibitem{SONG} X-C. Song,  Z. Phys.  {\bf C55}, 417 (1992)

\bibitem{OZ} O. Ogievetsky and B. Zumino, Lett. Math. Phys. {\bf 25}, 121
(1992)

\bibitem{KARPACZ} J.A. de Azc\'arraga and F. Rodenas, in {\it Quantum groups:
formalism and applications}, J. Lukierski, Z. Popowicz and J. Sobczyk eds.,
PWN (1994)

\bibitem{CL} A. Connes and J. Lott, Nucl. Phys. (Proc. Suppl.) {\bf B18},
29 (1990)

\bibitem{KS} D. Kastler, T. Sch\"ucker, {\it A detailed account of A. Connes'
version of the standard model IV}, CPT-94/P.3092 (Jan. 1995), hep-th 9501077

\end{thebibliography}
\end{document}